\documentclass[final,1p,times]{elsarticle}

\usepackage{amssymb}
 \usepackage{amsthm}
\usepackage{amscd}
\usepackage{amsmath}
\usepackage{amsfonts}
\usepackage{amssymb}
\usepackage{graphicx}
\newtheorem{theorem}{Theorem}
\newtheorem{proposition}[theorem]{Proposition}
\newtheorem{lemma}[theorem]{Lemma}

\usepackage{mathrsfs}
\usepackage{titletoc}

\newcommand{\la}{\Delta}

\newcommand{\ra}{\rightarrow}
\newcommand{\p}{\partial}
\newcommand{\f}{\frac}

\newcommand{\be}{\begin{equation}}
\renewcommand{\ra}{\rightarrow}
\newcommand{\ee}{\end{equation}}
\newcommand{\bea}{\begin{eqnarray}}
\newcommand{\eea}{\end{eqnarray}}
\newcommand{\bna}{\begin{eqnarray*}}
\newcommand{\ena}{\end{eqnarray*}}

\renewcommand{\le}{\left}
\newcommand{\ri}{\right}

\journal{***}

\begin{document}

\begin{frontmatter}

\title{Blow-up analysis concerning singular Trudinger-Moser inequalities in dimension two}

 \author{Yunyan Yang}
 \ead{yunyanyang@ruc.edu.cn}
 \author{Xiaobao Zhu}
 \ead{zhuxiaobao@ruc.edu.cn}

\address{ Department of Mathematics,
Renmin University of China, Beijing 100872, P. R. China}

\begin{abstract}
  In this paper, we derive a sharp version of the singular
  Trudinger-Moser inequality, which was originally established by Adimurthi and Sandeep
  (Nonlinear Differ. Equ. Appl. 2007). Moreover, extremal functions for those singular Trudinger-Moser inequalities
   are also obtained. Our method is the blow-up analysis. Compared with our previous work (J. Differential Equations 2015),
   the essential difficulty caused by the presence of singularity is how to analyse the asymptotic behaviour of
   certain maximizing sequence near the blow-up point. We overcome this difficulty by
  combining two different classification theorems of  Chen and  Li (Duke Math. J. 1991; Duke Math. J. 1995) to  get
  the desired bubble.
  \end{abstract}


\end{frontmatter}

\titlecontents{section}[0mm]
                       {\vspace{.2\baselineskip}}
                       {\thecontentslabel~\hspace{.5em}}
                        {}
                        {\dotfill\contentspage[{\makebox[0pt][r]{\thecontentspage}}]}
\titlecontents{subsection}[3mm]
                       {\vspace{.2\baselineskip}}
                       {\thecontentslabel~\hspace{.5em}}
                        {}
                       {\dotfill\contentspage[{\makebox[0pt][r]{\thecontentspage}}]}

\setcounter{tocdepth}{2}


\section{Introduction}

Let $\Omega$ be a  smooth bounded domain in $\mathbb{R}^2$, $W_0^{1,2}(\Omega)$ be a completion of
$C_0^\infty(\Omega)$ under the norm
$\|u\|_{W_0^{1,2}(\Omega)}=(\int_\Omega|\nabla u|^2dx)^{1/2}$.
 The Sobolev embedding theorem states that $W_0^{1,2}(\Omega)$ is embedded in $L^p(\Omega)$ for any $p>1$,
 but not in $L^\infty(\Omega)$.
However, as a limit case of the Sobolev embedding, the Trudinger-Moser inequality \cite{24,19,17,22,14} says that
\be\label{T-M}\sup_{u\in W_0^{1,2}(\Omega),\,\|\nabla u\|_2\leq 1}\int_\Omega e^{\gamma u^2}dx<+\infty,\quad \forall
\gamma\leq 4\pi;\ee
Moreover, these integrals are still finite for all $\gamma>4\pi$ and all $u\in W_0^{1,2}(\Omega)$, but the supremum is infinity.
This inequality was generalized in many ways, one of which is as below. Let $\lambda_1(\Omega)$
 be the first eigenvalue of the Laplacian, namely
 \be\label{first-eig}\lambda_1(\Omega)=\inf_{u\in W_0^{1,2}(\Omega),\,u\not\equiv 0}\f{\|\nabla u\|_2^2}{\|u\|_2^2}.\ee
 Here and throughout this paper, we denote the usual $L^p(\Omega)$-norm by $\|\cdot\|_p$ for any $p>0$.
  It was proved by Adimurthi and O. Druet \cite{A-D} that for any $\alpha<\lambda_1(\Omega)$,
  \be\label{a-d}\sup_{u\in W_0^{1,2}(\Omega),\,\|\nabla u\|_2\leq 1}\int_\Omega e^{4\pi u^2(1+\alpha\|u\|_2^2)}dx<+\infty;\ee
  Moreover the supremum is infinity for any $\alpha\geq \lambda_1(\Omega)$. This result was  extended by Y. Yang
  \cite{Yang-JFA,Yang-Tran} to the cases of
  high dimension and compact Riemannian surface, by Lu-Yang \cite{Lu-Yang} and J. Zhu \cite{ZhuJ}
  to the version of $L^p$-norm, by de Souza and J. M. do \'O \cite{de-O-1,doo1} to
  the whole Euclidean space, and by Tintarev \cite{Tint} to the following form
  \be\label{Tintarev}\sup_{u\in W_0^{1,2}(\Omega),\,\|u\|_{1,\alpha}\leq 1}\int_\Omega e^{4\pi u^2}dx<+\infty, \quad\forall
  \alpha<\lambda_1(\Omega).\ee
  Here and throughout this paper, for any $\alpha$ and $u$ satisfying $\sqrt{\alpha}\|u\|_2\leq \|\nabla u\|_2$, we denote
  \be\label{1-a}\|u\|_{1,\alpha}=\le(\int_\Omega|\nabla u|^2dx-\alpha\int_\Omega u^2dx\ri)^{1/2}.\ee
  One can check that (\ref{Tintarev}) is stronger that (\ref{a-d}). In a recent work \cite{Yang-JDE}, we
   generalized the inequality (\ref{Tintarev}) to the case that large eigenvalues are
  involved, as well as to the manifold case. Also, we obtained extremal functions for these kind of
  Trudinger-Moser inequalities. For pioneer works on extremal functions for Trudinger-Moser inequality,
  we refer the reader to L. Carleson and A. Chang \cite{CC}, M. Struwe \cite{Struwe},
  M. Flucher \cite{Flucher},  K. Lin \cite{Lin}, and Y. Li \cite{Lijpde}.\\

  Now we describe another kind of generalization of (\ref{T-M}), namely the singular Trudinger-Moser inequality.
  Based on a rearrangement argument,
  Adimurthi and K. Sandeep \cite{Adi-Sandeep} were able to prove the following:
  Let $\Omega\subset\mathbb{R}^2$ be a smooth bounded domain, and $0\leq \beta<1$ be fixed. Then there holds
  \be\label{Sing-T-M}\sup_{u\in W_0^{1,2}(\Omega),\,\|\nabla u\|_2\leq 1}\int_\Omega \f{e^{4\pi(1-\beta) u^2}}{|x|^{2\beta}}dx<+\infty.\ee
Clearly (\ref{Sing-T-M}) reduces to (\ref{T-M}) when $\beta=0$. This result was extended by Adimurthi and Y. Yang \cite{Adi-Yang}
to the whole Euclidean space, by de Souza and J. M. do \'O \cite{de-O-2} to another version in $\mathbb{R}^2$.
Such singular Trudinger-Moser inequalities are very important in the study of partial differential equations, see for examples
\cite{Adi-Sandeep,Adi-Yang,de-O-2,Yang-Jfa,Yang-jfa}.
When $\Omega$ is the unit ball $\mathbb{B}$, (\ref{Sing-T-M}) was improved by A. Yuan and X. Zhu \cite{YuanZhu} to the following:
Let $0\leq\beta<1$ be fixed, $\|u\|_{2,\beta}=(\int_\mathbb{B}|x|^{-2\beta}u^2dx)^{1/2}$, and
$$\lambda_{1,\beta}(\mathbb{B})=\inf_{u\in W_0^{1,2}(\mathbb{B}),\,\|u\|_{2,\beta}=1}\int_\mathbb{B}|\nabla u|^2dx.$$
Then we have for any $\alpha<\lambda_{1,\beta}(\mathbb{B})$,
\be\label{Yuan}\sup_{u\in W_0^{1,2}(\mathbb{B}),\,\|\nabla u\|_2\leq 1}\int_{\mathbb{B}}\f{e^{4\pi(1-\beta)u^2
(1+\alpha\|u\|_{2,\beta}^2)}}
{|x|^{2\beta}}dx<+\infty.\ee
Recently, an analog of (\ref{Yuan}) with $\|u\|_{2,\beta}$ replaced by $\|u\|_{p,\beta}$
was obtained by A. Yuan and Z. Huang \cite{YuanHuang}. The method of
\cite{YuanHuang,YuanZhu} is a symmetrization argument.\\

In this paper, we have two goals. One is to improve (\ref{Sing-T-M}) to a stronger version of
 the singular Trudinger-Moser inequality, namely, a combination of (\ref{Tintarev}) and (\ref{Sing-T-M}).
 Certainly, this improvement is also stronger than that of \cite{YuanHuang,YuanZhu} in our setting.
 The other is to prove the existence of extremal
functions for such stronger inequalities.
Our main results are stated as following:
\begin{theorem}\label{Thm1}
Let $\Omega\subset\mathbb{R}^2$ be a smooth bounded domain and $0\in\Omega$. Let $0<\beta<1$ be fixed
and  $\lambda_1(\Omega)$  be the first eigenvalue of the Laplacian with respect to Dirichlet boundary
condition given as in (\ref{first-eig}).
Then for any  $\alpha<\lambda_1(\Omega)$, the supremum
$$\sup_{u\in W_0^{1,2}(\Omega),\,\|u\|_{1,\alpha}\leq 1}\int_\Omega
\f{e^{4\pi(1-\beta) u^2}}{|x|^{2\beta}}dx<+\infty,$$
where $\|u\|_{1,\alpha}$ is defined as in (\ref{1-a}).
\end{theorem}

\begin{theorem}\label{Thm2}
Let $\Omega\subset\mathbb{R}^2$ be a smooth bounded domain and $0\in\Omega$. Let $0<\beta<1$ be fixed
and  $\lambda_1(\Omega)$  be the first eigenvalue of the Laplacian with respect to Dirichlet boundary
condition given as in (\ref{first-eig}).
Then for any  $\alpha<\lambda_1(\Omega)$ and any $\gamma\leq 4\pi(1-\beta)$, the supremum
$$\sup_{u\in W_0^{1,2}(\Omega),\,\|u\|_{1,\alpha}\leq 1}\int_\Omega
\f{e^{\gamma u^2}}{|x|^{2\beta}}dx$$
can be attained by some function $u_0\in W_0^{1,2}(\Omega)\cap C^1_{\rm loc}(\overline{\Omega}\setminus\{0\})\cap C^0(\overline\Omega)$
with $\|u_0\|_{1,\alpha}=1$,
where $\|u\|_{1,\alpha}$ is defined as in (\ref{1-a}).
\end{theorem}

The special case $\beta=0$ of Theorems \ref{Thm1} and \ref{Thm2} was already done by Y. Yang
via the method of blow-up analysis in \cite{Yang-JDE}.
   Though the only difference between \cite{Yang-JDE} and the current paper
   is the presence of the singular term $|x|^{-2\beta}$ with $0<\beta<1$,
   the previous blow-up procedure can not be applied directly. The essential difficulty caused by $|x|^{-2\beta}$
   is how to describe the exact asymptotic behavior of certain maximizing sequence near the blow-up point.
   Unlike in \cite{Yang-JDE}, we employ two different classification theorems
   of Chen and Li \cite{CL,CL2} to get the desired bubble. Of course, other steps of the
   blow-up analysis become more delicate because of the presence of $|x|^{-2\beta}$.
   The method of blow-up analysis is now a standard method
   of dealing with the best Trudinger-Moser inequalities. For works in this direction, we refer the reader to
    Carleson-Chang \cite{CC}, Struwe \cite{Struwe},
   Ding, Jost, Li and Wang \cite{DJLW}, Adimurthi and Struwe \cite{Adi-Stru}, Li \cite{Lijpde}, Adimurthi and Druet \cite{A-D}.
Using a concentration compactness alternative by Lions \cite{Lions}
 and following the lines of Flucher, Csato and Roy \cite{CR} proved the existence of extremal functions for the
  singular Trudinger-Moser inequality (\ref{Sing-T-M}).
   However, it seems that their method can not be applied to Theorems \ref{Thm1} and \ref{Thm2} for all range of $\alpha<\lambda_1(\Omega)$.

   Similarly as in \cite{Yang-JDE}, we may consider the case $\alpha\geq \lambda_1(\Omega)$. Note that the supremum in Theorem
    \ref{Thm1} is infinity in this case. Let $\lambda_1(\Omega)<\cdots<\lambda_j(\Omega)<\lambda_{j+1}(\Omega)<\cdots$ be all distinct eigenvalues of the Laplacian
   operator with respect to Dirichlet boundary condition. The corresponding eigenfunction space reads
$$E_{\lambda_j(\Omega)}=\le\{u\in W_0^{1,2}(\Omega): -\Delta u=\lambda_j(\Omega)u\ri\}.$$
Define $E_\ell=E_{\lambda_1(\Omega)}\oplus E_{\lambda_2(\Omega)}\oplus\cdots \oplus E_{\lambda_\ell(\Omega)}$ and
its orthogonal complement space in $W_0^{1,2}(\Omega)$ by
\be
\label{Eperp}E_\ell^\perp=\le\{u\in W_0^{1,2}(\Omega): \int_\Omega  u  vdx=0, \forall v\in E_\ell\ri\}.
\ee
It is known that for any positive integer $\ell$,
\be\label{1-1-1}\lambda_{\ell+1}(\Omega)=\inf_{u\in E_\ell^\perp,\,u\not\equiv 0}\f{\int_\Omega|\nabla u|^2dx}{\int_\Omega u^2dx},\ee
and that $\lambda_\ell\ra+\infty$ as $\ell\ra+\infty$.
For large $\alpha$, we have the following:

\begin{theorem}\label{Thm3}
Let $\Omega\subset\mathbb{R}^2$ be a smooth bounded domain and $0\in\Omega$. Let $0<\beta<1$ be fixed,
$\lambda_{\ell+1}(\Omega)$ be the $(\ell+1)$-th eigenvalue of the Laplacian operator given as in (\ref{1-1-1})
and $E_\ell^\perp$ be a function space defined as in (\ref{Eperp}).
Then for any  $\alpha<\lambda_{\ell+1}(\Omega)$, the supremum
$$\sup_{u\in E_\ell^\perp,\,\|u\|_{1,\alpha}\leq 1}\int_\Omega
\f{e^{4\pi(1-\beta) u^2}}{|x|^{2\beta}}dx<+\infty,$$
where $\|u\|_{1,\alpha}$ is defined as in (\ref{1-a}).
\end{theorem}

\begin{theorem}\label{Thm4}
Let $\Omega\subset\mathbb{R}^2$ be a smooth bounded domain and $0\in\Omega$. Let $0<\beta<1$ be fixed,
$\lambda_{\ell+1}(\Omega)$ be the $(\ell+1)$-th eigenvalue of the Laplacian operator given as in (\ref{1-1-1})
and $E_\ell^\perp$ be a function space defined as in (\ref{Eperp}).
Then for any  $\alpha<\lambda_{\ell+1}(\Omega)$ and any $\gamma\leq 4\pi(1-\beta)$, the supremum
$$\sup_{u\in  E_\ell^\perp,\,\|u\|_{1,\alpha}\leq 1}\int_\Omega
\f{e^{\gamma u^2}}{|x|^{2\beta}}dx$$
can be attained by some function $u_0\in  E_\ell^\perp\cap C^1_{\rm loc}(\overline{\Omega}\setminus\{0\})\cap C^0(\overline\Omega)$
with $\|u_0\|_{1,\alpha}=1$, where $\|u\|_{1,\alpha}$ is defined as in (\ref{1-a}).
\end{theorem}

The proof of Theorems \ref{Thm3} and \ref{Thm4} is completely analogous to that of Theorems \ref{Thm1} and \ref{Thm2}, except that
we must take effort to construct test functions $\phi_\epsilon\in E_\ell^\perp$ in the final step of the proof of Theorem \ref{Thm4}.
Our method of proving Theorems \ref{Thm1}-\ref{Thm4} can be applied to establish singular versions of other kind of Trudinger-Moser inequalities, say
the Hardy-Trudinger-Moser inequality \cite{Wang-Ye,Yang-Zhu-AGAG} and the Trudinger-Moser inequality
 involving the Gaussian curvature \cite{Yang-JGA}.

 We are informed by the referee that singular Trudinger-Moser inequalities for compact Riemannian surface
have been established by S. Iula and G. Mancini \cite{I-M} by using similar blow-up procedure. Also they obtained existence results of
extremal functions for those inequalities. It should be remarked that they derived
an upper bound of the singular Trudinger-Moser functional
by using Onofri's inequality (see \cite{I-M}, Theorem 1.1), while we deduced similar upper bound via the capacity
estimate (see Section \ref{2.3} below).  \\

The remaining part of this paper is organized as follows: In Section 2, we
use blow-up analysis to prove Theorems \ref{Thm1} and \ref{Thm2}; In Section 3, we prove Theorems \ref{Thm3}
and \ref{Thm4} by using a similar method. Throughout this paper, we do not distinguish sequence and subsequence.


\section{Proof of Theorems \ref{Thm1} and \ref{Thm2}}
We prove Theorems \ref{Thm1} and \ref{Thm2} jointly and divide the proof into several subsections.

\subsection{Maximizers for subcritical singular Trudinger-Moser functionals}

 We first show that maximizers for subcritical functionals exist. Namely,

\begin{proposition}\label{Prop1}
For any $\epsilon$, $0<\epsilon<1-\beta$, there exists some
$u_\epsilon\in W_0^{1,2}(\Omega)\cap C^1_{\rm loc}(\overline{\Omega}\setminus\{0\})\cap C^0(\overline\Omega)$ satisfying
$\|u_\epsilon\|_{1,\alpha}=1$ and
\be\label{subcritical}\int_\Omega
\f{e^{4\pi(1-\beta-\epsilon) u_\epsilon^2}}{|x|^{2\beta}}dx=\Lambda_{\beta,\epsilon}:=\sup_{u\in W_0^{1,2}(\Omega),\,\|u\|_{1,\alpha}\leq 1}\int_\Omega
\f{e^{4\pi(1-\beta-\epsilon) u^2}}{|x|^{2\beta}}dx.\ee
Moreover, in the distributional sense, $u_\epsilon$ satisfies the following equation
\be\label{EL}\le\{\begin{array}{lll}-\Delta u_\epsilon-\alpha u_\epsilon=\f{1}{\lambda_\epsilon}|x|^{-2\beta}u_\epsilon
e^{4\pi(1-\beta-\epsilon) u_\epsilon^2}\quad{\rm in}\quad \Omega,\\[1.2ex]
u_\epsilon\geq 0\quad {\rm in}\quad\Omega,
\\[1.2ex]
\lambda_\epsilon=\int_\Omega |x|^{-2\beta}u_\epsilon^2e^{4\pi(1-\beta-\epsilon) u_\epsilon^2}dx.
\end{array}\ri.\ee
\end{proposition}

\noindent{\it Proof.} Let $0<\epsilon<1-\beta$ be fixed. Take a function sequence $u_j\in W_0^{1,2}(\Omega)$ such that
$\|u_j\|_{1,\alpha}\leq 1$ and
\be\label{max}\int_\Omega
\f{e^{4\pi(1-\beta-\epsilon) u_j^2}}{|x|^{2\beta}}dx\ra \sup_{u\in W_0^{1,2}(\Omega),\,\|u\|_{1,\alpha}\leq 1}\int_\Omega
\f{e^{4\pi(1-\beta-\epsilon) u^2}}{|x|^{2\beta}}dx\ee
as $j\ra\infty$.
Since $\alpha<\lambda_1(\Omega)$, we have
$$\le(1-\f{\alpha}{\lambda_1(\Omega)}\ri)\int_\Omega|\nabla u_j|^2dx\leq \int_\Omega|\nabla u_j|^2dx-\alpha\int_\Omega u_j^2dx\leq 1,$$
 which implies that $u_j$ is bounded in $W_0^{1,2}(\Omega)$. Hence there exists some $u_\epsilon\in W_0^{1,2}(\Omega)$
such that up to a subsequence, we have $u_j\rightharpoonup u_\epsilon$ weakly in $W_0^{1,2}(\Omega)$, $u_j\ra u_\epsilon$
strongly in $L^q(\Omega)$ for any $q>1$, and $u_j\ra u_\epsilon$ almost everywhere  in $\Omega$.
For any $1<p<1/\beta$, $\delta>0$, $s>1$ and $s^\prime=s/(s-1)$, we have by the H\"older inequality,
\bea\nonumber
\int_\Omega |x|^{-2\beta p}e^{4\pi(1-\beta-\epsilon)p u_j^2}dx&\leq&\int_\Omega
|x|^{-2\beta p}e^{4\pi(1-\beta-\epsilon)p(1+\delta)(u_j-u_\epsilon)^2+4\pi(1-\beta-\epsilon)p(1+\f{1}{4\delta})u_\epsilon^2}dx\\
\nonumber&\leq&\le(\int_\Omega|x|^{-2\beta p}e^{4\pi(1-\beta-\epsilon)p(1+\delta)s(u_j-u_\epsilon)^2}dx\ri)^{1/s}\\\label{lp}
&&\times\le(\int_\Omega|x|^{-2\beta p}e^{4\pi(1-\beta-\epsilon)p(1+\f{1}{4\delta})s^\prime u_\epsilon^2}dx\ri)^{1/{s^\prime}}.
\eea
Choosing $p$, $1+\delta$ and $s$ sufficiently close to $1$, we have
\be\label{less}(1-\beta-\epsilon)p(1+\delta)s+\beta p<1.\ee
Note that
\bea\nonumber
\int_\Omega|\nabla (u_j-u_\epsilon)|^2dx&=&\int_\Omega|\nabla u_j|^2dx-\int_\Omega|\nabla u_\epsilon|^2dx+o_j(1)\\
\nonumber&=&\|u_j\|_{1,\alpha}^2-\|u_\epsilon\|_{1,\alpha}^2+o_j(1)\\
&\leq&1-\|u_\epsilon\|_{1,\alpha}^2+o_j(1),\label{energ}
\eea
since $\|u_j\|_{1,\alpha}\leq 1$. Inserting (\ref{less}) and (\ref{energ}) into (\ref{lp}), we have by the singular
Trudinger-Moser inequality (\ref{Sing-T-M}) that $|x|^{-2\beta}e^{4\pi(1-\beta-\epsilon)u_j^2}$ is bounded in $L^p(\Omega)$
for some $p>1$. Since
\bna
|x|^{-2\beta}|e^{4\pi(1-\beta-\epsilon) u_j^2}-e^{4\pi(1-\beta-\epsilon) u_\epsilon^2}|\leq
4\pi(1-\beta-\epsilon)|x|^{-2\beta}(e^{4\pi(1-\beta-\epsilon) u_j^2}+e^{4\pi(1-\beta-\epsilon) u_\epsilon^2})|u_j^2-u_\epsilon^2|
\ena
and $u_j\ra u_\epsilon$ strongly in $L^q(\Omega)$ for all $q>1$ as $j\ra\infty$, we conclude that
\be\label{lim}\lim_{j\ra+\infty}\int_\Omega|x|^{-2\beta}e^{4\pi(1-\beta-\epsilon)u_j^2}dx=
\int_\Omega|x|^{-2\beta}e^{4\pi(1-\beta-\epsilon)u_\epsilon^2}dx.\ee
It follows from (\ref{energ}) that
\be\label{e-1}\|u_\epsilon\|_{1,\alpha}\leq 1.\ee
Combining (\ref{max}), (\ref{lim}) and (\ref{e-1}), we have that $u_\epsilon$ attains the supremum $\Lambda_{\beta,\epsilon}$.
Clearly $u_\epsilon\not\equiv 0$.
Suppose $\|u_\epsilon\|_{1,\alpha}<1$. It follows that
$$\Lambda_{\beta,\epsilon}
=\int_\Omega |x|^{-2\beta}e^{4\pi(1-\beta-\epsilon)u_\epsilon^2}dx<
\int_\Omega |x|^{-2\beta}e^{4\pi(1-\beta-\epsilon)\f{u_\epsilon^2}{\|u_\epsilon\|_{1,\alpha}^2}}dx\leq \Lambda_{\beta,\epsilon},$$
which is a contradiction. Hence
we have $\|u_\epsilon\|_{1,\alpha}=1$. Also one can see that $|u_\epsilon|$ attains the supremum $\Lambda_{\beta,\epsilon}$. Hence
$u_\epsilon$ can be chosen such that $u_\epsilon\geq 0$. A straightforward calculation shows that $u_\epsilon$ satisfies the Euler-Lagrange
equation (\ref{EL}). $\hfill\Box$\\

In view of Proposition \ref{Prop1}, to prove Theorem \ref{Thm2}, we only need to prove that there exists some function
$u^\ast\in W_0^{1,2}(\Omega)\cap C^1_{\rm loc}(\overline{\Omega}\setminus\{0\})\cap C^0(\overline{\Omega})$ verifying that
$u^\ast\geq 0$, $\|u^\ast\|_{1,\alpha}=1$, and
\be\label{extremal}\int_\Omega |x|^{-2\beta}e^{4\pi(1-\beta){u^\ast}^2}dx=\sup_{u\in W_0^{1,2}(\Omega),\,\|u\|_{1,\alpha}\leq 1}\int_\Omega
|x|^{-2\beta}e^{4\pi(1-\beta)u^2}dx.\ee

\subsection{Blow-up analysis}

Since $u_\epsilon$ is bounded in $W_0^{1,2}(\Omega)$, we can assume without loss of generality,
\bea\label{ueps-weak-con}
&&u_\epsilon\rightharpoonup u_0\quad{\rm weakly\,\,in}\quad W_0^{1,2}(\Omega),\\[1.2ex]\label{strong}
&&u_\epsilon\ra u_0\quad{\rm strongly\,\,in}\quad L^q(\Omega),\,\,\,\forall q\geq 1,
\\[1.2ex]\label{ae-con}
&&u_\epsilon\ra u_0\quad{\rm a.\,e.\,\,\,in}\quad\Omega.
\eea
Let $c_\epsilon=\max_{\Omega}u_\epsilon$. If $c_\epsilon$ is bounded, then for any $u\in W_0^{1,2}(\Omega)$ with
$\|u\|_{1,\alpha}\leq 1$, we have by the Lebesgue dominated convergence theorem
\bna\int_\Omega |x|^{-2\beta}e^{4\pi(1-\beta) u^2}dx&=&\lim_{\epsilon\ra 0}
\int_\Omega |x|^{-2\beta}e^{4\pi(1-\beta-\epsilon)u^2}dx\\
&\leq&\lim_{\epsilon\ra 0}
\int_\Omega |x|^{-2\beta}e^{4\pi(1-\beta-\epsilon)u_\epsilon^2}dx\\
&=&\int_\Omega |x|^{-2\beta}e^{4\pi(1-\beta)u_0^2}dx.\ena
Hence $u_0$ is the desired maximizer, or equivalently (\ref{extremal}) holds. In the following, we can assume
$c_\epsilon=u_\epsilon(x_\epsilon)\ra +\infty$ and $x_\epsilon\ra x_0\in\overline{\Omega}$ as $\epsilon\ra 0$.
By an inequality $e^{t^2}\leq 1+t^2e^{t^2}$, we have
\bna
\int_\Omega|x|^{-2\beta} e^{4\pi(1-\beta-\epsilon)u_\epsilon^2}dx\leq
\int_\Omega|x|^{-2\beta}dx+4\pi\lambda_\epsilon.
\ena
This together with (\ref{subcritical}) leads to
\be\label{lower}\liminf_{\epsilon\ra 0}\lambda_\epsilon>0.\ee
\begin{proposition}\label{Lionslemma} We have
$u_0\equiv0$, $x_0=0$, and $|\nabla u_\epsilon|^2dx\rightharpoonup \delta_0$, where $\delta_0$ denotes the Dirac
measure centered at $0$.
\end{proposition}
\noindent{\it Proof.} Suppose $u_0\not\equiv 0$, then we have $$\int_\Omega|\nabla(u_\epsilon-u_0)|^2dx
=1-\|u_0\|_{1,\alpha}^2+o_\epsilon(1).$$
In view of (\ref{lower}) and a similar estimate as (\ref{lp}), we have by applying elliptic estimates to (\ref{EL}), $u_\epsilon$ is bounded
in $W^{2,p}(\Omega)$ for some $p>1$. Hence the Sobolev embedding theorem implies that $u_\epsilon$ is bounded in
$C^0(\overline{\Omega})$. In particular, $c_\epsilon$ is bounded, contradicting $c_\epsilon\ra+\infty$ as $\epsilon\ra 0$.
Hence $u_0\equiv 0$.

Since $\int_\Omega|\nabla u_\epsilon|^2dx=1+o_\epsilon(1)$, it is not difficult to see that
$|\nabla u_\epsilon|^2dx\rightharpoonup \delta_{x_0}$, for otherwise
we have by using elliptic estimates, $u_\epsilon$ is uniformly bounded near $x_0$.
This contradicts again $c_\epsilon\ra+\infty$ as $\epsilon\ra 0$. Moreover, we have $u_\epsilon\ra 0$ in
$C^1_{\rm loc}(\overline{\Omega}\setminus\{0,\,x_0\})\cap C^0_{\rm loc}(\overline{\Omega}\setminus\{x_0\})$.

Suppose $x_0\not=0$. Then $\lambda_\epsilon^{-1}|x|^{-2\beta}u_\epsilon e^{4\pi(1-\beta-\epsilon)u_\epsilon^2}$ is
bounded in $L^{q_1}(\mathbb{B}_{|x_0|/2})$ for some $q_1>1$. Noting that $|x|^{-2\beta}\leq (|x_0|/2)^{-2\beta}$ when $|x|\geq |x_0|/2$,
we have that $\lambda_\epsilon^{-1}|x|^{-2\beta}u_\epsilon e^{4\pi(1-\beta-\epsilon)u_\epsilon^2}$ is
bounded in $L^{q_2}(\Omega\setminus\mathbb{B}_{|x_0|/2})$ for some $q_2>1$. Therefore $\lambda_\epsilon^{-1}|x|^{-2\beta}u_\epsilon e^{4\pi(1-\beta-\epsilon)u_\epsilon^2}$ is
bounded in $L^q(\Omega)$ for $q=\min\{q_1,q_2\}>1$.
 Hence we have by using elliptic estimates, $c_\epsilon$ is bounded contradicting
$c_\epsilon\ra+\infty$. This completes the proof of the proposition. $\hfill\Box$\\

Let
\be\label{scal}r_\epsilon=\sqrt{\lambda_\epsilon}c_\epsilon^{-1}e^{-2\pi(1-\beta-\epsilon)c_\epsilon^2}.\ee
Note that $u_\epsilon\ra 0$ in $L^q(\Omega)$ for any $q\geq 1$. We obtain
\be\label{tends-0}r_\epsilon^2 e^{4\pi\tau c_\epsilon^2}\ra 0,\quad \forall \tau<1-\beta.\ee
 We now distinguish two cases to proceed.\\

{\it Case 1.} $|x_\epsilon|^{1-\beta}/r_\epsilon\ra +\infty$.

Define on $\Omega_{1,\epsilon}=\{x\in\mathbb{R}^2: x_\epsilon+r_\epsilon|x_\epsilon|^\beta x\in\Omega\}$,
$$w_\epsilon(x)=c_\epsilon^{-1}u_\epsilon(x_\epsilon+r_\epsilon|x_\epsilon|^\beta x),\quad
v_\epsilon(x)=c_\epsilon(u_\epsilon(x_\epsilon+r_\epsilon|x_\epsilon|^\beta x)-c_\epsilon).$$
A straightforward calculation shows
\be\label{lap-w}-\Delta w_\epsilon(x)=\alpha r_\epsilon^2|x_\epsilon|^{2\beta}w_\epsilon+c_\epsilon^{-2}|x_\epsilon|^{2\beta}
\le|x_\epsilon+r_\epsilon|x_\epsilon|^\beta x\ri|^{-2\beta}w_\epsilon e^{4\pi(1-\beta-\epsilon)c_\epsilon^2(w_\epsilon^2-1)}\quad
{\rm in}\quad \Omega_{1,\epsilon}.\ee
Since $0\leq w_\epsilon\leq 1$ and $|x_\epsilon|^{2\beta}
|x_\epsilon+r_\epsilon|x_\epsilon|^\beta x|^{-2\beta}=1+o_\epsilon(1)$, where $o_\epsilon(1)\ra 0$ in $\mathbb{B}_R$
for any $R>0$,
we have by applying elliptic estimates to (\ref{lap-w}) that $w_\epsilon\ra w$ in $C^1_{\rm loc}(\mathbb{R}^2)$, where $w$ satisfies
$$-\Delta w(x)=0\quad{\rm in}\quad \mathbb{R}^2.$$
Since $w\leq 1$ and $w(0)=1$, the Liouville theorem leads to $w\equiv 1$. Also we have
\be\label{lap-v}-\Delta v_\epsilon=\alpha c_\epsilon^2r_\epsilon^2|x_\epsilon|^{2\beta}w_\epsilon+|x_\epsilon|^{2\beta}
\le|x_\epsilon+r_\epsilon|x_\epsilon|^\beta x\ri|^{-2\beta}w_\epsilon e^{4\pi(1-\beta-\epsilon)(w_\epsilon+1)v_\epsilon}\quad{\rm in}\quad \Omega_{1,\epsilon}.\ee
Clearly we have by applying elliptic estimates to (\ref{lap-v}) that $v_\epsilon\ra v$ in $C^1_{\rm loc}(\mathbb{R}^2)$, where $v$ satisfies
\be\label{bubbel}
   \le\{\begin{array}{lll}&-\la
     v=e^{8\pi(1-\beta)v}\quad{\rm in}\quad\mathbb{R}^2,\\[1.2 ex]
     &v(0)=0=\sup_{\mathbb{R}^2}v.
     \end{array}
    \ri.\ee
On one hand, we have for any $R>0$,
\bna\nonumber
\int_{\mathbb{B}_R(0)}e^{8\pi(1-\beta)v}dx&=&\lim_{\epsilon\ra 0}\int_{\mathbb{B}_R(0)}
e^{4\pi(1-\beta-\epsilon)(u_\epsilon^2(x_\epsilon+r_\epsilon|x_\epsilon|^\beta x)-c_\epsilon^2)}dx\\\nonumber
&=&\lim_{\epsilon\ra 0}\lambda_\epsilon^{-1}\int_{\mathbb{B}_{Rr_\epsilon|x_\epsilon|^\beta}(x_\epsilon)}|x_\epsilon|^{-2\beta}
c_\epsilon^2e^{4\pi(1-\beta-\epsilon)u_\epsilon^2(y)}dy\\\nonumber
&=&\lim_{\epsilon\ra 0}\lambda_\epsilon^{-1}\int_{\mathbb{B}_{Rr_\epsilon|x_\epsilon|^\beta}(x_\epsilon)}|y|^{-2\beta}
u_\epsilon^2(y)e^{4\pi(1-\beta-\epsilon)u_\epsilon^2(y)}dy\\
&\leq&1.
\ena
This leads to
\be\label{leq-1} \int_{\mathbb{R}^2}e^{8\pi(1-\beta)v}dx\leq 1.\ee
On the other hand, in view of (\ref{bubbel}) and (\ref{leq-1}), a result of Chen and Li \cite{CL} implies that
$v$ is radially symmetric and
\be\label{geq-1}\int_{\mathbb{R}^2}e^{8\pi (1-\beta)v}dx\geq \f{1}{1-\beta}.\ee
The contradiction between (\ref{leq-1}) and (\ref{geq-1}) indicates that
{\it Case 1} can not occur.\\

{\it Case 2.} $|x_\epsilon|^{1-\beta}/r_\epsilon\leq C$ for some constant $C$.

Denote $\Omega_\epsilon=\{x\in\mathbb{R}^2: x_\epsilon+r_\epsilon^{{1}/{(1-\beta)}}x\in\Omega\}$. Define
$$\psi_\epsilon(x)=c_\epsilon^{-1}u_\epsilon(x_\epsilon+r_\epsilon^{1/(1-\beta)} x),\quad
\varphi_\epsilon(x)=c_\epsilon(u_\epsilon(x_\epsilon+r_\epsilon^{1/(1-\beta)} x)-c_\epsilon).$$
It follows that $\psi_\epsilon$ is a distributional solution to the equation
\be\label{lap-psi}-\Delta\psi_\epsilon=\alpha r_\epsilon^{\f{2}{1-\beta}}\psi_\epsilon+c_\epsilon^{-2} |x+r_\epsilon^{-1/(1-\beta)}x_\epsilon|^{-2\beta}\psi_\epsilon
e^{4\pi(1-\beta-\epsilon)(1+\psi_\epsilon)\varphi_\epsilon}\quad{\rm in}\quad \Omega_\epsilon.\ee
 In view of (\ref{tends-0}), $r_\epsilon\ra 0$ and thus
$\Omega_\epsilon\ra\mathbb{R}^2$. We can assume $r_\epsilon^{-1/(1-\beta)}x_\epsilon\ra x^\ast$ for some $x^\ast\in\mathbb{R}^2$.
Applying elliptic estimates to (\ref{lap-psi}), we have that $\psi_\epsilon\ra \psi_0$ in $C^1_{\rm loc}(\mathbb{R}^2\setminus\{-x^\ast\})
\cap C^0_{\rm loc}(\mathbb{R}^2)$,
where $\psi_0$ is a distributional harmonic function on $\mathbb{R}^2$. Since
$\psi_0(x)\leq\limsup_{\epsilon\ra 0}\psi_\epsilon(x)\leq 1$ for all $x\in\mathbb{R}^2$ and $\psi_0(0)=\lim_{\epsilon\ra 0}\psi_\epsilon(0)=1$,
 the Liouville theorem implies that $\psi_0\equiv 1$
on $\mathbb{R}^2$. Hence we conclude
\be\label{ten-1}\psi_\epsilon\ra 1\quad{\rm in}\quad C^1_{\rm loc}(\mathbb{R}^2\setminus\{-x^\ast\})\cap C^0_{\rm loc}(\mathbb{R}^2).\ee
Clearly, $\varphi_\epsilon$ is a distributional solution to
\be\label{lap-varphi}-\Delta\varphi_\epsilon=\alpha c_\epsilon^2r_\epsilon^{\f{2}{1-\beta}}\psi_\epsilon+|x+r_\epsilon^{-1/(1-\beta)}x_\epsilon|^{-2\beta}
\psi_\epsilon
e^{4\pi(1-\beta-\epsilon)(1+\psi_\epsilon)\varphi_\epsilon}\quad{\rm in}\quad \Omega_\epsilon.\ee
Since $\psi_\epsilon\ra 1$ in $C^1_{\rm loc}(\mathbb{R}^2\setminus\{-x^\ast\})$ and
$\varphi_\epsilon(0)=0=\max_{\mathbb{R}^2}\varphi_\epsilon$, applying elliptic estimates to (\ref{lap-varphi}), we have that
$\varphi_\epsilon\ra \varphi_0$ in $C^1_{\rm loc}(\mathbb{R}^2\setminus\{-x^\ast\})\cap C^0_{\rm loc}(\mathbb{R}^2)$, where $\varphi_0$ is a solution to
\be\label{varphi}-\Delta \varphi_0=|x+x^\ast|^{-2\beta}e^{8\pi(1-\beta)\varphi_0}\quad{\rm in}\quad \mathbb{R}^2\setminus\{-x^\ast\}.\ee
If we let $y=x_\epsilon+r_\epsilon^{1/(1-\beta)}x$ with $|x+x^\ast|\leq R$, then for any fixed $R>|x^\ast|+1$,
there holds $|y|\leq 2Rr_\epsilon^{1/(1-\beta)}$. Combining (\ref{ten-1}) and Fatou's lemma, we have
\bna
 \int_{\mathbb{B}_R(-x^\ast)}|x+x^\ast|^{-2\beta}e^{8\pi(1-\beta)\varphi_0}dx&\leq&
 \limsup_{\epsilon\ra 0}\int_{\mathbb{B}_R(-x^\ast)}|x+r_\epsilon^{-1/(1-\beta)}x_\epsilon|^{-2\beta}e^{4\pi(1-\beta-\epsilon)
 (1+\psi_\epsilon)\varphi_\epsilon}dx\\
 &\leq&\limsup_{\epsilon\ra 0}\lambda_\epsilon^{-1}\int_{\mathbb{B}_{2Rr_\epsilon^{1/(1-\beta)}}(0)}
 |y|^{-2\beta}u_\epsilon^2(y)e^{4\pi(1-\beta-\epsilon) u_\epsilon^2(y)}dy\\
 &\leq&1.
\ena
Hence
$$\int_{\mathbb{R}^2}|x+x^\ast|^{-2\beta}e^{8\pi(1-\beta)\varphi_0}dx\leq 1.$$
 By a classification result of Chen and Li (\cite{CL2}, Theorem 3.1), we have
\be\label{decphi}\varphi_0(x)=-\f{1}{4\pi(1-\beta)}\log
\le(1+\f{\pi}{1-\beta}|x+x^\ast|^{2(1-\beta)}\ri).\ee
Note that
\be\label{decphi1}\varphi_0(0)=\lim_{\epsilon\ra 0}\varphi_\epsilon(0)=0.\ee
Combining (\ref{decphi}) and (\ref{decphi1}), we have that $x^\ast=0$ and thus
\be\label{dec-phi1}\varphi_0(x)=-\f{1}{4\pi(1-\beta)}\log
\le(1+\f{\pi}{1-\beta}|x|^{2(1-\beta)}\ri).\ee
It follows that
\be\label{integ-1}\int_{\mathbb{R}^2}|x|^{-2\beta}e^{8\pi(1-\beta)\varphi_0}dx=1.\ee

Define $u_{\epsilon,\gamma}=\min\{u_\epsilon,\gamma c_\epsilon\}$. Similar to \cite{Lijpde,A-D}, we have the following:
\begin{lemma}\label{le-1}
For any $\gamma$, $0<\gamma<1$, there holds
$$\lim_{\epsilon\ra 0}\int_\Omega|\nabla u_{\epsilon,\gamma}|^2dx=\gamma.$$
\end{lemma}

\noindent{\it Proof.} In view of the equation (\ref{EL}), we have by using the integration by parts,
  \bna
  \int_\Omega|\nabla u_{\epsilon,\gamma}|^2dx&=&\int_\Omega\nabla u_{\epsilon,\gamma}\nabla u_\epsilon dx
  =-\int_\Omega u_{\epsilon,\gamma}\la u_\epsilon
  dx\\ &=&\lambda_\epsilon^{-1}\int_\Omega |x|^{-2\beta}u_{\epsilon,\gamma}
  u_\epsilon e^{4\pi(1-\beta-\epsilon)u_\epsilon^2} dx+\alpha\int_\Omega u_\epsilon u_{\epsilon,\gamma}dx\\
  &\geq&\lambda_\epsilon^{-1}\int_{\mathbb{B}_{Rr_\epsilon^{{1}/{(1-\beta)}}}(x_\epsilon)}|x|^{-2\beta}u_{\epsilon,\gamma}
  u_\epsilon e^{4\pi(1-\beta-\epsilon)u_\epsilon^2} dx+o_\epsilon(1)\\
  &=&\gamma(1+o_\epsilon(1))\int_{\mathbb{B}_R(0)}|y+r_\epsilon^{-1/(1-\beta)}x_\epsilon|^{-2\beta}e^{4\pi(1-\beta-\epsilon)
  (u_\epsilon^2(x_\epsilon+r_\epsilon^{1/(1-\beta)}y)-c_\epsilon^2)}dy+o_\epsilon(1),
  \ena
  which leads to
  $$\liminf_{\epsilon\ra 0}\int_\Omega|\nabla u_{\epsilon,\gamma}|^2dx\geq \gamma
  \int_{\mathbb{B}_R(0)}|y|^{-2\beta}e^{8\pi(1-\beta)\varphi(y)}dy, \quad\forall R>0.$$
  In view of (\ref{integ-1}), we have by passing to the limit $R\ra +\infty$ in the above inequality,
  \be\label{gee}\liminf_{\epsilon\ra 0}\int_\Omega|\nabla u_{\epsilon,\gamma}|^2dx\geq\gamma.\ee
  Note that $|\nabla(u_\epsilon-\gamma c_\epsilon)^+|^2=
  \nabla(u_\epsilon-\gamma c_\epsilon)^+\nabla u_\epsilon$ on $\Omega$ and
  $(u_\epsilon-\gamma c_\epsilon)^+=(1+o_\epsilon(1))(1-\gamma)c_\epsilon$
  on $\mathbb{B}_{Rr_\epsilon^{1/(1-\beta)}}(x_\epsilon)$.  Similarly as above, we obtain
  \be\label{ggg}\liminf_{\epsilon\ra 0}\int_\Omega|\nabla (u_\epsilon-\gamma c_\epsilon)^+|^2dx\geq 1-\gamma.\ee
  Since $|\nabla u_\epsilon|^2=|\nabla u_{\epsilon,\gamma}|^2+|\nabla(u_\epsilon-\gamma c_\epsilon)^+|^2$ almost
  everywhere, we get
  \be\label{111}\int_\Omega |\nabla u_{\epsilon,\gamma}|^2dx+
  \int_\Omega|\nabla (u_\epsilon-\gamma c_\epsilon)^+|^2dx=\|u_\epsilon\|_{1,\alpha}^2+\alpha\int_\Omega u_\epsilon^2dx=1+o_\epsilon(1).
  \ee
  Combining (\ref{gee}), (\ref{ggg}), and (\ref{111}), we finish the proof of the lemma. $\hfill\Box$\\

  As a consequence of Lemma \ref{le-1}, we have the following:
  \begin{lemma}\label{u-p}
  There holds
  $$\lim_{\epsilon\ra 0}\int_\Omega|x|^{-2\beta}e^{4\pi(1-\beta-\epsilon)u_\epsilon^2}dx\leq
  \int_\Omega|x|^{-2\beta}dx+\limsup_{\epsilon\ra 0}\f{\lambda_\epsilon}{c_\epsilon^2}.$$
  \end{lemma}
  \noindent{\it Proof.} For any $\gamma$, $0<\gamma<1$, there holds
  \bea\nonumber
  \int_\Omega|x|^{-2\beta}e^{4\pi(1-\beta-\epsilon)u_\epsilon^2}dx&=&
  \int_{u_\epsilon\leq \gamma c_\epsilon}|x|^{-2\beta}e^{4\pi(1-\beta-\epsilon)u_\epsilon^2}dx+
  \int_{u_\epsilon>\gamma c_\epsilon}|x|^{-2\beta}e^{4\pi(1-\beta-\epsilon)u_\epsilon^2}dx\\\label{u1}
  &\leq&\int_\Omega|x|^{-2\beta}e^{4\pi(1-\beta-\epsilon)u_{\epsilon,\gamma}^2}dx
  +\f{\lambda_\epsilon}{\gamma^2c_\epsilon^2}.
  \eea
  By Lemma \ref{le-1}, $|x|^{-2\beta}e^{4\pi(1-\beta-\epsilon)u_{\epsilon,\gamma}^2}$ is bounded in $L^q(\Omega)$ for
  some $q>1$. Note also that $u_{\epsilon,\gamma}$ converges to $0$ almost everywhere. Hence
  $|x|^{-2\beta}e^{4\pi(1-\beta-\epsilon)u_{\epsilon,\gamma}^2}$ converges to $|x|^{-2\beta}$ in $L^1(\Omega)$. Passing to the
  limit $\epsilon\ra 0$ in (\ref{u1}), we obtain
  $$\lim_{\epsilon\ra 0}\int_\Omega|x|^{-2\beta}e^{4\pi(1-\beta-\epsilon)u_\epsilon^2}dx\leq
  \int_\Omega|x|^{-2\beta}dx+\f{1}{\gamma^2}\limsup_{\epsilon\ra 0}\f{\lambda_\epsilon}{c_\epsilon^2}.$$
  Letting $\gamma\ra 1$, we conclude the lemma. $\hfill\Box$\\

  It follows from Lemma \ref{u-p} that
  \be\label{t-0}\limsup_{\epsilon\ra 0}\f{\lambda_\epsilon}{c_\epsilon^\theta}=+\infty,\quad\forall\theta<2.\ee
  For otherwise, we have $\lambda_\epsilon/c_\epsilon^2\ra 0$ as $\epsilon\ra 0$. Let $v\in W_0^{1,2}(\Omega)$ be such that
  $\|v\|_{1,\alpha}=1$. Then we have by Lemma \ref{u-p}  that
  \bna\int_\Omega|x|^{-2\beta}e^{4\pi(1-\beta) v^2}dx&\leq&\sup_{u\in W_0^{1,2}(\Omega),\,\|u\|_{1,\alpha}\leq 1}\int_\Omega
|x|^{-2\beta}e^{4\pi(1-\beta) u^2}dx\\&=&
\lim_{\epsilon\ra 0}\int_\Omega|x|^{-2\beta}e^{4\pi(1-\beta-\epsilon)u_\epsilon^2}dx\\
&=&\int_\Omega |x|^{-2\beta}dx.\ena
This is impossible since $v\not\equiv 0$. Thus (\ref{t-0}) holds.\\

\begin{lemma}\label{delta-0} $c_\epsilon u_\epsilon$ is bounded in $W_0^{1,q}(\Omega)$ for any $1<q<2$.
Furthermore, $c_\epsilon u_\epsilon\rightharpoonup G$ weakly in $W_0^{1,q}(\Omega)$ for any $1<q<2$ and
$c_\epsilon u_\epsilon\ra G$ strongly in $L^r(\Omega)$ for any $r>1$, where $G$ satisfies
\be\label{green}-\Delta G-\alpha G=\delta_0\ee
in the distributional sense, $\delta_0$ stands for the Dirac measure centered at $0$.
\end{lemma}
\noindent{\it Proof.} Firstly we claim that for any $\phi\in C^2(\overline{\Omega})$, there holds
\be\label{d-0}\lim_{\epsilon\ra 0}
\int_\Omega \lambda_\epsilon^{-1}|x|^{-2\beta}c_\epsilon u_\epsilon e^{4\pi(1-\beta-\epsilon) u_\epsilon^2}\phi dx
=\phi(0).\ee
To see this, we denote $g_\epsilon=\lambda_\epsilon^{-1}|x|^{-2\beta}c_\epsilon u_\epsilon e^{4\pi(1-\beta-\epsilon) u_\epsilon^2}$.
Clearly
\be\label{intest}\int_\Omega g_\epsilon\phi dx=\int_{u_\epsilon<\gamma c_\epsilon}g_\epsilon\phi dx+
\int_{\{u_\epsilon\geq\gamma c_\epsilon\}\setminus\mathbb{B}_{Rr_\epsilon^{1/(1-\beta)}}(x_\epsilon)}g_\epsilon\phi dx
+\int_{\mathbb{B}_{Rr_\epsilon^{1/(1-\beta)}}(x_\epsilon)\cap\{u_\epsilon\geq\gamma c_\epsilon\}}g_\epsilon\phi dx.\ee
We estimate the three integrals on the right hand of (\ref{intest}) respectively. By (\ref{t-0}) and Lemma \ref{le-1},
\be\label{I1}\int_{u_\epsilon<\gamma c_\epsilon}g_\epsilon\phi dx=\f{c_\epsilon}{\lambda_\epsilon}\int_{u_\epsilon<\gamma c_\epsilon}
|x|^{-2\beta}u_\epsilon e^{4\pi(1-\beta-\epsilon) u_{\epsilon,\gamma}^2}\phi dx=o_\epsilon(1).
\ee
Since $\mathbb{B}_{Rr_\epsilon^{1/(1-\beta)}}(x_\epsilon)\subset\{u_\epsilon\geq \gamma c_\epsilon\}$ for sufficiently small
$\epsilon>0$, we have by (\ref{integ-1}),
\bea\nonumber
\int_{\mathbb{B}_{Rr_\epsilon^{1/(1-\beta)}}(x_\epsilon)\cap\{u_\epsilon\geq\gamma c_\epsilon\}}g_\epsilon\phi dx&=&
\phi(0)(1+o_\epsilon(1))\int_{\mathbb{B}_{Rr_\epsilon^{1/(1-\beta)}}(x_\epsilon)}\lambda_\epsilon^{-1} c_\epsilon|x|^{-2\beta}
u_\epsilon e^{4\pi(1-\beta-\epsilon) u_\epsilon^2}dx\\\nonumber
&=&\phi(0)(1+o_\epsilon(1))\int_{\mathbb{B}_R(0)}|x|^{-2\beta}e^{8\pi(1-\beta)\varphi}dx\\\label{I2}
&=&\phi(0)(1+o_\epsilon(1)+o_R(1)).
\eea
Noting that
\bna
\int_{\{u_\epsilon\geq\gamma c_\epsilon\}\setminus\mathbb{B}_{Rr_\epsilon^{1/(1-\beta)}}(x_\epsilon)}g_\epsilon dx
&\leq& \f{1}{\gamma}\int_{\{u_\epsilon\geq\gamma c_\epsilon\}\setminus\mathbb{B}_{Rr_\epsilon^{1/(1-\beta)}}(x_\epsilon)}\lambda_\epsilon^{-1}|x|^{-2\beta}u_\epsilon^2
e^{4\pi(1-\beta-\epsilon) u_\epsilon^2}dx\\
&\leq&\f{1}{\gamma}\le(1-\int_{\mathbb{B}_{Rr_\epsilon^{1/(1-\beta)}}(x_\epsilon)}\lambda_\epsilon^{-1}|x|^{-2\beta}u_\epsilon^2
e^{4\pi(1-\beta-\epsilon) u_\epsilon^2}dx\ri)\\
&=&\f{1}{\gamma}\le(1-\int_{\mathbb{B}_R(0)}|x|^{-2\beta}e^{8\pi(1-\beta)\varphi}dx\ri),
\ena
we have
\be\label{I3}\lim_{R\ra+\infty}\lim_{\epsilon\ra 0}\int_{\{u_\epsilon\geq\gamma c_\epsilon\}\setminus
\mathbb{B}_{Rr_\epsilon^{1/(1-\beta)}}(x_\epsilon)}g_\epsilon\phi dx=0.\ee
Inserting (\ref{I1})-(\ref{I3}) to (\ref{intest}), we conclude (\ref{d-0}).

By the equation (\ref{EL}), $c_\epsilon u_\epsilon$ is a distributional solution to
\be\label{g-epsi}-\Delta(c_\epsilon u_\epsilon)-\alpha c_\epsilon u_\epsilon=g_\epsilon\quad{\rm in}\quad\Omega.\ee
It follows from (\ref{d-0}) that $g_\epsilon$ is bounded in $L^1(\Omega)$.
We {\it claim} that $c_\epsilon u_\epsilon$ is bounded in $L^1(\Omega)$. To see this, we suppose on the contrary,
$\|c_\epsilon u_\epsilon\|_{L^1(\Omega)}\ra +\infty$ as $\epsilon\ra 0$.
Define a new sequence of functions $\chi_\epsilon=c_\epsilon u_\epsilon/\|c_\epsilon u_\epsilon\|_{L^1(\Omega)}$.
Then applying a result of Struwe (\cite{Str-1}, Theorem 2.2) to (\ref{g-epsi}), we have that $\chi_\epsilon$ is bounded in $W_0^{1,q}(\Omega)$ for any $q$, $1<q<2$,
in particular $\chi_\epsilon\ra \chi$ strongly in $L^1(\Omega)$. Since $g_\epsilon/\|c_\epsilon u_\epsilon\|_{L^1(\Omega)}\ra 0$
in $L^1(\Omega)$, $\chi$ is a distributional solution to
$-\Delta \chi-\alpha\chi=0$ in $\Omega$, which leads to $\chi\equiv 0$. This contradicts $\|\chi\|_{L^1(\Omega)}=\lim_{\epsilon\ra 0}
\|\chi_\epsilon\|_{L^1(\Omega)}=1$ and confirms our claim. Now since $g_\epsilon+\alpha c_\epsilon u_\epsilon$ is bounded in $L^1(\Omega)$,
applying again (\cite{Str-1}, Theorem 2.2) to (\ref{g-epsi}), we conclude that
 $c_\epsilon u_\epsilon$ is bounded in $W_0^{1,q}(\Omega)$ for any $q$, $1<q<2$.
Hence there exists some $G\in \cap_{1<q<2}W_0^{1,q}(\Omega)$ such that $c_\epsilon u_\epsilon\rightharpoonup G$
weakly in $W_0^{1,q}(\Omega)$ for any $1<q<2$, and that $c_\epsilon u_\epsilon\ra G$ strongly in $L^r(\Omega)$ for any $r>1$.
Since (\ref{d-0}) implies that $g_\epsilon\rightharpoonup \delta_0$ in sense of meaure, where $\delta_0$ denotes the Dirac measure centered at
$0$. In view of (\ref{g-epsi}), $G$ is a distributional solution to (\ref{green}). $\hfill\Box$\\

Obviously, $G$ takes the form
\be\label{gr-dec}G(x)=-\f{1}{2\pi}\log|x|+A_0+\psi(x),\ee
where $A_0$ is a constant and $\psi\in C^1(\overline{\Omega})$.

\subsection{An upper bound}\label{2.3}

In this subsection, we use the capacity estimate, which was first used by Y. Li \cite{Lijpde} in this topic,  to derive an upper bound
of the integrals  $\int_\Omega |x|^{-2\beta}e^{4\pi(1-\beta-\epsilon)u_\epsilon^2}dx$.
 Take small $\delta$ such that $B_{2\delta}(0)\subset \Omega$.
 Define a function space $$\mathscr{W}_\epsilon(a,b)=\{u\in W^{1,2}(\mathbb{B}_\delta(x_\epsilon)\setminus \mathbb{B}
 _{Rr_\epsilon^{1/(1-\beta)}}(x_\epsilon)):
 u|_{\p \mathbb{B}_\delta(x_\epsilon)}=a,\, u|_{\p
 \mathbb{B}_{Rr_\epsilon^{1/(1-\beta)}}(x_\epsilon)}=b\}.$$
 It is not difficult to see that
 $\inf\limits_{u\in \mathscr{W}_\epsilon(s_\epsilon,i_\epsilon)}\int_{\mathbb{B}_\delta(x_\epsilon)\setminus
 \mathbb{B}_{Rr_\epsilon^{1/(1-\beta)}}(x_\epsilon)}|\nabla u|^2dx$
 is attained by $h(x)$ satisfying
 $$\le\{\begin{array}{llll}&\Delta h=0\quad{\rm in} \quad
     \mathbb{B}_\delta(x_\epsilon)\setminus \mathbb{B}_{Rr_\epsilon^{1/(1-\beta)}}(x_\epsilon)\\
     &h|_{\p \mathbb{B}_\delta(x_\epsilon)}=s_\epsilon\\
     &h|_{\p \mathbb{B}_{Rr_\epsilon^{1/(1-\beta)}}(x_\epsilon)}=i_\epsilon.
     \end{array}
     \ri.$$
     One can check that
     $$h(x)=\f{s_\epsilon(\log|x-x_\epsilon|-\log(Rr_\epsilon^{1/(1-\beta)}))+i_\epsilon
     (\log\delta-\log|x-x_\epsilon|)}{\log\delta-\log(Rr_\epsilon^{1/(1-\beta)})},$$
     and that
     \be\label{3.1}\int_{\mathbb{B}_\delta(x_\epsilon)\setminus \mathbb{B}_{Rr_\epsilon^{1/(1-\beta)}}(x_\epsilon)}|\nabla
     h|^2dx
     =\f{2\pi(s_\epsilon-i_\epsilon)^2}{\log\delta-\log(Rr_\epsilon^{1/(1-\beta)})}.\ee
     Let
     $$i_\epsilon=\inf_{\p \mathbb{B}_{Rr_\epsilon^{1/(1-\beta)}}(x_\epsilon)}u_\epsilon,\quad s_\epsilon=\sup_{\p\mathbb{B}_\delta(x_\epsilon)}
     u_\epsilon$$
     and
     $\widetilde{u}_\epsilon=\max\{s_\epsilon,\min\{u_\epsilon,i_\epsilon\}\}$.
     Then $\widetilde{u}_\epsilon\in \mathscr{W}_\epsilon(s_\epsilon,i_\epsilon)$ and
     $|\nabla \widetilde{u}_\epsilon|\leq |\nabla u_\epsilon|$ a. e. in
     $\mathbb{B}_\delta(x_\epsilon)\setminus \mathbb{B}_{Rr_\epsilon^{1/(1-\beta)}}(x_\epsilon)$,
     provided that $\epsilon$ is chosen sufficiently small. Hence
     \bea
     \int_{\mathbb{B}_\delta(x_\epsilon)\setminus \mathbb{B}_{Rr_\epsilon^{1/(1-\beta)}}(x_\epsilon)}|\nabla h|^2
     dx&\leq&
     \int_{\mathbb{B}_\delta(x_\epsilon)\setminus \mathbb{B}_{Rr_\epsilon^{1/(1-\beta)}}(x_\epsilon)}
     |\nabla \widetilde{u}_\epsilon|^2dx
     \leq\int_{\mathbb{B}_\delta(x_\epsilon)\setminus \mathbb{B}_{Rr_\epsilon^{1/(1-\beta)}}(x_\epsilon)}
     |\nabla u_\epsilon|^2dx\nonumber\\[1.2ex]
     &\leq& 1+\alpha\int_\Omega u_\epsilon^2dx-\int_{\Omega\setminus \mathbb{B}_\delta(x_\epsilon)}|\nabla
     u_\epsilon|^2dx
     -\int_{\mathbb{B}_{Rr_\epsilon^{1/(1-\beta)}}(x_\epsilon)}|\nabla
     u_\epsilon|^2dx.\quad
     \eea
     Now we compute $\int_{\mathbb{B}_{Rr_\epsilon}(x_\epsilon)}|\nabla
     u_\epsilon|^2dx$ and $\int_{\Omega\setminus \mathbb{B}_\delta(x_\epsilon)}
     |\nabla u_\epsilon|^2dx$. In view of (\ref{gr-dec}) and (\ref{green}),
     integration by parts leads to
     \bna
     \int_{\Omega\setminus \mathbb{B}_\delta(x_\epsilon)}|\nabla G|^2dx&=&-\int_{\Omega\setminus B_\delta(x_\epsilon)}G\Delta Gdx
     -\int_{\p\mathbb{B}_{\delta}(x_\epsilon)}G\f{\p G}{\p\nu}ds\\
     &=&-\f{1}{2\pi}\log\delta+A_0+\alpha\int_\Omega G^2dx+o_\epsilon(1)+o_\delta(1).
     \ena
  Since $c_\epsilon u_\epsilon\ra G$ in $C^1_{\rm loc}(\overline{\Omega}\setminus\{0\})$,
we obtain
\be\label{out}\int_{\Omega\setminus \mathbb{B}_\delta(x_\epsilon)}|\nabla u_\epsilon|^2dx=
\f{1}{c_\epsilon^2}\le(-\f{1}{2\pi}\log\delta+A_0+\alpha\int_\Omega G^2dx+
o_\epsilon(1)+o_\delta(1)\ri).\ee
Let $\varphi_0$ be given as in (\ref{dec-phi1}). A straightforward calculation shows
$$\int_{\mathbb{B}_R(0)}|\nabla\varphi_0|^2dx=\f{1}{2\pi}\log R+\f{1}{4\pi(1-\beta)}
\log\f{\pi}{1-\beta}-\f{1}{4\pi(1-\beta)}+O(\f{1}{R^{2-2\beta}}).$$
 Since $\varphi_\epsilon\ra \varphi_0$ in $W^{2,q}_{\rm loc}(\mathbb{R}^2)$ for some $q>1$, in particular in $W^{1,2}_{\rm loc}
 (\mathbb{R}^2)$, we have
 \bna\int_{\mathbb{B}_{Rr_\epsilon^{1/(1-\beta)}}(x_\epsilon)}|\nabla u_\epsilon|^2dx&=&
  \int_{\mathbb{B}_R(0)}c_\epsilon^{-2}|\nabla\varphi_\epsilon(y)|^2dy\\
  &=&\f{1}{c_\epsilon^2}
\le(\int_{\mathbb{B}_R(0)}|\nabla\varphi_0(y)|^2dy+o_\epsilon(1)\ri)\\
  &=&
 \f{1}{c_\epsilon^2}
\le(\f{1}{2\pi}\log R+\f{1}{4\pi(1-\beta)}
\log\f{\pi}{1-\beta}-\f{1}{4\pi(1-\beta)}+o(1)\ri),\ena
where $o(1)\ra 0$ as $\epsilon\ra 0$ first and then $R\ra+\infty$ .
It follows from (\ref{dec-phi1}) and Lemma \ref{delta-0} that
 \bna
 i_\epsilon&=&c_\epsilon+\f{1}{c_\epsilon}\le(-\f{1}{4\pi(1-\beta)}\log\le(1+\f{\pi}{1-\beta}
 R^{2(1-\beta)}\ri)+o(1)\ri),\\[1.2ex]
 s_\epsilon&=&\f{1}{c_\epsilon}\le(-\f{1}{2\pi}\log\delta+A_0+o(1)\ri),
\ena
where $o(1)\ra 0$ as $\epsilon\ra 0$ first and then $\delta\ra 0$. Hence
 $$2\pi(s_\epsilon-i_\epsilon)^2=2\pi c_\epsilon^2-\f{1}{1-\beta}\log
 \le(1+\f{\pi}{1-\beta}R^{2(1-\beta)}\ri)+2\log\delta-4\pi A_0+o(1).$$
 Recalling (\ref{scal}), we have
\be\label{3.5}\log\delta-\log(Rr_\epsilon^{1/(1-\beta)})=\log\delta-\log R
-\f{1}{2(1-\beta)}\log\f{\lambda_\epsilon}{c_\epsilon^2}+\f{2\pi(1-\beta-\epsilon)
c_\epsilon^2}{(1-\beta)}.\ee Combining
(\ref{3.1})--(\ref{3.5}) and noting that
$$\alpha\int_\Omega u_\epsilon^2dx=\f{1}{c_\epsilon}\le(\int_\Omega G^2dx+o_\epsilon(1)\ri),$$
we have
\bna
&&2\pi c_\epsilon^2-\f{1}{1-\beta}\log
 \le(1+\f{\pi}{1-\beta}R^{2(1-\beta)}\ri)+2\log\delta-4\pi A_0+o(1)\\
 &&\leq \le(\log\delta-\log R
-\f{1}{2(1-\beta)}\log\f{\lambda_\epsilon}{c_\epsilon^2}+\f{2\pi(1-\beta-\epsilon)
c_\epsilon^2}{(1-\beta)}+o_\epsilon(1)\ri)\\
&&\quad\times\le(1-\f{1}{c_\epsilon^2}\le(\f{\log R}{2\pi}+\f{\log\f{\pi}{1-\beta}}{4\pi(1-\beta)}-
\f{1}{4\pi(1-\beta)}-\f{\log\delta}{2\pi}+A_0+o(1)\ri)\ri)\\
&&=\f{2\pi(1-\beta-\epsilon)
c_\epsilon^2}{(1-\beta)}+2\log\delta-2\log R-\f{1+o(1)}{2(1-\beta)}\log\f{\lambda_\epsilon}{c_\epsilon^2}
-\f{\log\f{\pi}{1-\beta}}{2(1-\beta)}\\
&&\quad+\f{1}{2(1-\beta)}-2\pi A_0+o(1).\ena
It follows that
$$\f{1+o(1)}{2(1-\beta)}\log\f{\lambda_\epsilon}{c_\epsilon^2}\leq \f{\log\f{\pi}{1-\beta}}{2(1-\beta)}+
\f{1}{2(1-\beta)}+2\pi A_0+o(1),$$
which implies that
$$\limsup_{\epsilon\ra 0}\f{\lambda_\epsilon}{c_\epsilon^2}\leq \f{\pi}{1-\beta}e^{1+4\pi(1-\beta)A_0}.$$
Therefore we conclude by Lemma \ref{u-p},

\be\label{upper-bd}\sup_{u\in W_0^{1,2}(\Omega),\,\|u\|_{1,\alpha}\leq 1}\int_\Omega \f{e^{4\pi(1-\beta)u^2}}{|x|^{2\beta}}dx
=\lim_{\epsilon\ra 0}\int_\Omega\f{e^{4\pi(1-\beta-\epsilon)u_\epsilon^2}}{|x|^{2\beta}}dx\leq
\int_\Omega|x|^{-2\beta}dx+\f{\pi}{1-\beta}e^{1+4\pi(1-\beta)A_0}.\ee

\subsection{Completion of the proof of Theorem \ref{Thm1}}

Let $u_0$ be as in (\ref{ueps-weak-con})-(\ref{ae-con}). In case $c_\epsilon\ra+\infty$,  (\ref{upper-bd}) holds.
In case $c_\epsilon$ is bounded,  $u_0$ satisfies (\ref{extremal}).  In conclusion, there necessarily holds
$$\sup_{u\in W_0^{1,2}(\Omega),\,\|u\|_{1,\alpha}\leq 1}\int_\Omega
\f{e^{4\pi(1-\beta)u^2}}{|x|^{2\beta}}dx<+\infty.$$
This completes the proof of Theorem \ref{Thm1}. $\hfill\Box$

\subsection{Completion of the proof of Theorem \ref{Thm2}}

   In view of Proposition \ref{Prop1}, to finish the proof of Theorem \ref{Thm2}, we only need to prove (\ref{extremal}).
   If $c_\epsilon$ is bounded, then (\ref{extremal}) is already true. If $c_\epsilon\ra+\infty$, then (\ref{upper-bd}) holds.
   We shall construct a sequence of functions $\phi_\epsilon\in W_0^{1,2}(\Omega)$ with $\|\phi_\epsilon\|_{1,\alpha}=1$ such that
   $$\int_\Omega \f{e^{4\pi(1-\beta)\phi_\epsilon^2}}{|x|^{2\beta}}dx>\int_\Omega|x|^{-2\beta}dx+\f{\pi}{1-\beta}e^{1+4\pi(1-\beta)A_0}.$$
   This contradicts (\ref{upper-bd}). Hence $c_\epsilon$ must be bounded and the proof of Theorem \ref{Thm2} is finished.

     Define a sequence of functions on $\Omega$ by
     \be\label{ppp}\phi_\epsilon(x)=\le\{
     \begin{array}{llll}
     &c+\f{1}{c}\le(-\f{1}{4\pi(1-\beta)}\log(1+\f{\pi}{1-\beta}\f{|x|^{2(1-\beta)}}
     {\epsilon^{2(1-\beta)}})+b\ri),
     \quad &x\in\overline{\mathbb{B}}_{R\epsilon}\\[1.5ex]
     &\f{G-\eta \psi}{c},\quad & x\in \mathbb{B}_{2R\epsilon}\setminus \overline{\mathbb{B}}_{R\epsilon}\\[1.2ex]
     &\f{G}{c},\quad & x\in \Omega\setminus\mathbb{B}_{2R\epsilon},
     \end{array}
     \ri.\ee
     where $G$ and $\psi$ are functions given as in (\ref{gr-dec}), $R=(-\log\epsilon)^{1/(1-\beta)}$,
     $\eta\in C_0^1(\mathbb{B}_{2R\epsilon})$ satisfying that $\eta=1$ on $\mathbb{B}_{R\epsilon}$ and
     $|\nabla \eta|\leq \f{2}{R\epsilon}$, $b$ and $c$ are constants depending only on $\epsilon$ to be
     determined later. Here and in the sequel, $\mathbb{B}_r$ stands for a ball centered at $0$ with radius $r$.
     Clearly $\mathbb{B}_{2R\epsilon}\subset\Omega$ provided that $\epsilon$ is sufficiently small.
     In order to assure that $\phi_\epsilon\in W_0^{1,2}(\Omega)$, we set
     $$
     c+\f{1}{c}\le(-\f{1}{4\pi(1-\beta)}\log(1+\f{\pi}{1-\beta} R^{2(1-\beta)})+b\ri)
     =\f{1}{c}\le(-\f{1}{2\pi}\log (R\epsilon)+A_{0}\ri),
     $$
     which gives
     \be\label{2pic2-1}
      c^2=-\f{1}{2\pi}\log\epsilon+A_0-b+\f{1}{4\pi(1-\beta)}\log \f{\pi}{1-\beta}
     +O(\f{1}{R^{2(1-\beta)}}).
     \ee
     Noting that $\psi(x)=O(|x|)$ as $x\ra 0$, we have
     $|\nabla (\eta\psi)|=O(1)$ as $\epsilon\ra 0$. It follows that
     $$\int_{\mathbb{B}_{2R\epsilon}\setminus\mathbb{B}_{R\epsilon}}|\nabla(\eta\psi)|^2dx=O(R^2\epsilon^2),\quad
     \int_{\mathbb{B}_{2R\epsilon}\setminus\mathbb{B}_{R\epsilon}}\nabla G\nabla(\eta\psi)dx=O(R\epsilon).$$
     Integration by parts gives
     \bna
     \int_{\Omega\setminus\mathbb{B}_{R\epsilon}}|\nabla G|^2dx&=&-\int_{\Omega\setminus\mathbb{B}_{R\epsilon}}
     G\Delta Gdx-\int_{\p\mathbb{B}_{R\epsilon}}G\f{\p G}{\p\nu}ds\\
     &=&-\f{1}{2\pi}\log(R\epsilon)+A_0+\alpha\int_\Omega G^2dx+O(R\epsilon).
     \ena
     This leads to
     \bna
     \int_{\Omega\setminus\mathbb{B}_{R\epsilon}}|\nabla\phi_\epsilon|^2dx&=&\f{1}{c^2}\int_{\Omega\setminus\mathbb{B}_{R\epsilon}}
     |\nabla G|^2dx+\f{1}{c^2}\int_{\mathbb{B}_{2R\epsilon}\setminus\mathbb{B}_{R\epsilon}}|\nabla(\eta\psi)|^2dx\\
     &&
     -\f{2}{c^2}\int_{\mathbb{B}_{2R\epsilon}\setminus\mathbb{B}_{R\epsilon}}\nabla G\nabla(\eta\psi)dx\\
     &=&\f{1}{c^2}\le(-\f{1}{2\pi}\log(R\epsilon)+A_0+\alpha\int_\Omega G^2dx+O(R\epsilon)\ri).
     \ena
     Also we have
     \bna
     \int_{\mathbb{B}_{R\epsilon}}|\nabla\phi_\epsilon|^2dx&=&\f{1}{4(1-\beta)^2c^2}\int_{\mathbb{B}_{R}}
     \f{|x|^{2-4\beta}}{(1+\f{\pi}{1-\beta}|x|^{2(1-\beta)})^2}dx\\
     &=&\f{\pi}{2(1-\beta)^2c^2}\int_0^{R}\f{r^{3-4\beta}}{(1+\f{\pi}{1-\beta}r^{2-2\beta})^2}dr\\
     &=&\f{1}{4\pi(1-\beta)c^2}\int_0^{\f{\pi}{1-\beta}R^{2-2\beta}}\f{tdt}{(1+t)^2}\\
     &=&\f{1}{4\pi(1-\beta)c^2}\le(\log\f{\pi}{1-\beta}-1+\log R^{2-2\beta}+O(\f{1}{R^{2-2\beta}})\ri).
     \ena
     Hence
     $$\int_\Omega |\nabla\phi_\epsilon|^2dx=\f{1}{c^2}
     \le(-\f{\log\epsilon}{2\pi}+A_0+\alpha\int_\Omega G^2dx-\f{1}{4\pi(1-\beta)}+\f{1}{4\pi(1-\beta)}\log\f{\pi}{1-\beta}+O(\f{1}{R^{2-2\beta}})\ri).$$
     Note that
     $$\int_\Omega \phi_\epsilon^2dx=\f{1}{c^2}\le(\int_\Omega G^2dx+O(R\epsilon)\ri).$$
     Set
     $$\|\phi_\epsilon\|_{1,\alpha}^2=\int_\Omega|\nabla\phi_\epsilon|^2dx-\alpha\int_\Omega\phi_\epsilon^2dx=1.$$
     It follows that
     \be\label{c2-2}c^2=-\f{1}{2\pi}\log\epsilon+A_0-\f{1}{4\pi(1-\beta)}
     +\f{1}{4\pi(1-\beta)}\log\f{\pi}{1-\beta}+O(\f{1}{R^{2-2\beta}}).\ee
     Combining (\ref{2pic2-1}) and (\ref{c2-2}), we obtain
     \be\label{B}b=\f{1}{4\pi(1-\beta)}+O(\f{1}{R^{2-2\beta}}).\ee
     In view of (\ref{c2-2}) and (\ref{B}), there holds on $\mathbb{B}_{R\epsilon}$,
     \bna
     4\pi(1-\beta)\phi_\epsilon^2&\geq& 4\pi(1-\beta)c^2-2\log\le(1+\f{\pi}{1-\beta}\f{|x|^{2(1-\beta)}}{\epsilon^{2(1-\beta)}}\ri)
     +8\pi(1-\beta)b\\
     &=&-2\log\le(1+\f{\pi}{1-\beta}\f{|x|^{2(1-\beta)}}{\epsilon^{2(1-\beta)}}\ri)
     -2(1-\beta)\log\epsilon+1\\&&+4\pi(1-\beta)A_0
     +\log\f{\pi}{1-\beta}+O(\f{1}{R^{2-2\beta}}),
     \ena
     which together with the estimate
     \bna
     \int_{\mathbb{B}_{R}}\f{1}{(1+\f{\pi}{1-\beta}|y|^{2(1-\beta)})^2|y|^{2\beta}}dy&=&
     \int_0^R\f{2\pi r^{1-2\beta}}{(1+\f{\pi}{1-\beta}r^{2-2\beta})^2}dr\\
     &=&\int_0^{\f{\pi}{1-\beta}R^{2-2\beta}}\f{dt}{(1+t)^2}\\
     &=&1-\f{1}{1+\f{\pi}{1-\beta}R^{2-2\beta}}
     \ena
     leads to
     \bna
     \int_{\mathbb{B}_{R\epsilon}}|x|^{-2\beta}e^{4\pi(1-\beta)\phi_\epsilon^2}dx&\geq&
     \f{\pi}{1-\beta}\epsilon^{-2(1-\beta)}e^{1+4\pi(1-\beta)A_0+O(\f{1}{R^{2-2\beta}})}\\
     &&\times\int_{\mathbb{B}_{R\epsilon}}\f{1}{(1+\f{\pi}{1-\beta}\f{|x|^{2(1-\beta)}}{\epsilon^{2(1-\beta)}})^2|x|^{2\beta}}dx\\
     &=&\f{\pi}{1-\beta}e^{1+4\pi(1-\beta)A_0+O(\f{1}{R^{2-2\beta}})}\\
     &&\times\int_{\mathbb{B}_{R}}\f{1}{(1+\f{\pi}{1-\beta}|y|^{2(1-\beta)})^2|y|^{2\beta}}dy\\
     &=&\f{\pi}{1-\beta}e^{1+4\pi(1-\beta)A_0}+O(\f{1}{R^{2-2\beta}}).
     \ena

     On the other hand, since
     $$\int_{\mathbb{B}_{2R\epsilon}}|x|^{-2\beta}dx=O\le((R\epsilon)^{2-2\beta}\ri)=O(\f{1}{R^{2-2\beta}})$$
     and
     $$\int_{\mathbb{B}_{2R\epsilon}}|x|^{-2\beta}G^2dx=O\le((R\epsilon)^{2-2\beta}\log^2(R\epsilon)\ri)
     =O(\f{1}{R^{2-2\beta}}),$$
     we obtain
     \bna
     \int_{\Omega\setminus\mathbb{B}_{R\epsilon}}|x|^{-2\beta}e^{4\pi(1-\beta)\phi_\epsilon^2}dx&\geq&
     \int_{\Omega\setminus\mathbb{B}_{2R\epsilon}}|x|^{-2\beta}(1+4\pi(1-\beta)\phi_\epsilon^2)dx\\
     &=&\int_{\Omega\setminus\mathbb{B}_{2R\epsilon}}|x|^{-2\beta}dx+\f{4\pi(1-\beta)}{c^2}
     \int_{\Omega\setminus\mathbb{B}_{2R\epsilon}}|x|^{-2\beta}G^2dx\\
     &=&\int_\Omega|x|^{-2\beta}dx+\f{4\pi(1-\beta)}{c^2}
     \int_{\Omega}|x|^{-2\beta}G^2dx+O(\f{1}{R^{2-2\beta}}).
     \ena
     Therefore
     $$\int_{\Omega}|x|^{-2\beta}e^{4\pi(1-\beta)\phi_\epsilon^2}dx\geq\int_\Omega|x|^{-2\beta}dx+\f{\pi}{1-\beta}e^{1+4\pi(1-\beta)A_0}
     +\f{4\pi(1-\beta)}{c^2}
     \int_{\Omega}|x|^{-2\beta}G^2dx+O(\f{1}{R^{2-2\beta}}).$$
     In view of $R=(-\log\epsilon)^{1/(1-\beta)}$ and (\ref{c2-2}), we have
     $\f{1}{R^{2-2\beta}}=o(\f{1}{c^2})$, and thus
     $$\int_{\Omega}|x|^{-2\beta}e^{4\pi(1-\beta)\phi_\epsilon^2}dx>\int_\Omega|x|^{-2\beta}dx+\f{\pi}{1-\beta}e^{1+4\pi(1-\beta)A_0},$$
     provided that $\epsilon>0$ is chosen sufficiently small.
     $\hfill\Box$

\section{Proof of Theorems \ref{Thm3} and \ref{Thm4}}

Since the proof of Theorems \ref{Thm3} and \ref{Thm4} is analogous to that of Theorems \ref{Thm1} and \ref{Thm2},
we only give its outline, but emphasize the difference between them as below.

\subsection{Proof of Theorem \ref{Thm3}}

Since $E_\ell$ is a finite dimensional linear space, there exists an orthogonal basis
$\{\psi_1,\cdots,\psi_m\}$, namely $E_\ell={\rm Span}\{\psi_1,\cdots,\psi_m\}$, where
$\psi_i\in C^1(\overline{\Omega})$ satisfies
\be\label{orth}\int_\Omega\psi_i\psi_jdx=\delta_{ij}=\le\{\begin{array}{lll}
1,\,\,i=j\\[1.2ex]
0,\,\,i\not=j.
\end{array}\ri.\ee
Let $E_\ell^\perp$ be defined as in (\ref{Eperp}). Then it follows that
 $$E_\ell^\perp=\le\{u\in W_0^{1,2}(\Omega): \int_\Omega  u  \psi_i dx=0, 1\leq i\leq m\ri\}.$$
 Clearly $W_0^{1,2}(\Omega)=E_\ell\oplus E_\ell^\perp$. Note that $E_\ell^\perp$ is weakly closed, namely, if $u_k\in E_\ell^\perp$
 and $u_k\rightharpoonup u$ weakly in $W_0^{1,2}(\Omega)$, then $u\in E_\ell^\perp$.
     Using the argument of the proof of Proposition \ref{Prop1}, we can show that
     for any $\epsilon$, $0<\epsilon<1-\beta$, there exists
     some $u_\epsilon\in E_\ell^\perp$ with
 $\|u_\epsilon\|_{1,\alpha}=1$ such that
 \be\label{subcr-1}\int_\Omega \f{e^{4\pi(1-\beta-\epsilon)u_\epsilon^2}}{|x|^{2\beta}}dx=
 \sup_{u\in E_\ell^\perp,\,\|u\|_{1,\alpha}\leq 1}\int_\Omega \f{e^{4\pi(1-\beta-\epsilon)u^2}}
 {|x|^{2\beta}}dx.\ee Moreover $u_\epsilon$ satisfies the
 Euler-Lagrange equation
  $$\label{E-L-1-2}\le\{
  \begin{array}{lll}
  -\Delta u_\epsilon-\alpha u_\epsilon=\f{1}{\lambda_\epsilon}|x|^{-2\beta}u_\epsilon e^{(4\pi-\epsilon)u_\epsilon^2}-
  \sum_{i=1}^m\f{\gamma_\epsilon^i}{\lambda_\epsilon}\psi_i\,\,\,{\rm in}\,\,\,
  \Omega,\\[1.5ex] u_\epsilon\in E_\ell^\perp\cap C^1_{\rm loc}(\overline{\Omega}\setminus\{0\})\cap C^0(\overline{\Omega}),\\[1.5ex]
  \lambda_\epsilon=\int_\Omega |x|^{-2\beta}u_\epsilon^2 e^{(4\pi-\epsilon)u_\epsilon^2}dx,\\[1.5ex]
  \gamma_\epsilon^i=\int_\Omega \psi_i|x|^{-2\beta}u_\epsilon e^{(4\pi-\epsilon)u_\epsilon^2}dx.
  \end{array}
  \ri.$$
  Without loss of generality we can assume
  $u_\epsilon\rightharpoonup u_0$ weakly in $W_0^{1,2}(\Omega)$, $u_\epsilon\ra u_0$ strongly in $L^p(\Omega)$, $\forall p>1$,
  and $u_\epsilon\ra u_0$  a.\,e. in $\Omega$. Clearly we have $u_0\in E_\ell^\perp$ and $\|u_0\|_{1,\alpha}\leq 1$.
  If $u_\epsilon$ is bounded in $C^0(\overline\Omega)$, then for all $u\in E_\ell^\perp$ with
  $\|u\|_{1,\alpha}\leq 1$, we have by (\ref{subcr-1}) and the Lebesgue dominated convergence theorem
  $$\int_\Omega \f{e^{4\pi(1-\beta) {u}^2}}{|x|^{2\beta}}dx=\lim_{\epsilon\ra 0}\int_\Omega \f{e^{4\pi(1-\beta-\epsilon){u}^2}}
  {|x|^{2\beta}}dx
  \leq\lim_{\epsilon\ra 0}\int_{\Omega}\f{e^{4\pi(1-\beta-\epsilon){u_\epsilon}^2}}{|x|^{2\beta}}dx=\int_\Omega
  \f{e^{4\pi(1-\beta) u_0^2}}{|x|^{2\beta}}dx.$$
  Hence we have
  $$\label{b-d}\int_\Omega \f{e^{4\pi(1-\beta) u_0^2}}{|x|^{2\beta}}dx=
  \sup_{u\in E_\ell^\perp,\,\|u\|_{1,\alpha}\leq 1}\int_\Omega \f{e^{4\pi(1-\beta) u^2}}{|x|^{2\beta}}dx.$$
  It is easy to see that $\|u_0\|_{1,\alpha}=1$.
  The regularity theory implies that $u_0\in C_{\rm loc}^1(\overline{\Omega}\setminus\{0\})\cap C^0(\overline{\Omega})$,
  and thus $u_0$ is a desired extremal function.
  In the sequel we assume up to a subsequence
  $$c_\epsilon=\max_{\overline\Omega}|u_\epsilon|=\|u_\epsilon\|_{C^0(\overline\Omega)}
  \ra +\infty\quad{\rm as}\quad \epsilon\ra 0.$$
  Without loss of generality
  we assume $c_\epsilon=u_\epsilon(x_\epsilon)$. For otherwise we replace $u_\epsilon$ by $-u_\epsilon$ below.
  Then up to a subsequence, we can easily see that $x_\epsilon\ra 0$, $u_0\equiv 0$ and
  $|\nabla u_\epsilon|^2dx\rightharpoonup \delta_{0}$ weakly in sense of measure.
  Define a sequence of blow-up functions
  $$\varphi_\epsilon(x)=c_\epsilon(u_\epsilon(x_\epsilon+r_\epsilon x)-c_\epsilon)$$ for $x\in\Omega_\epsilon
  =\{x\in \mathbb{R}^2:x_\epsilon+r_\epsilon x\in\Omega\}$, where
 $r_\epsilon=\sqrt{\lambda_\epsilon}c_\epsilon^{-1}e^{-2\pi(1-\beta-\epsilon)c_\epsilon^2}$
 is the blow-up scale. Then up to a subsequence, there holds $\varphi_\epsilon\ra \varphi_0$  in
 $C^1_{\rm loc}(\mathbb{R}^2\setminus\{0\})\cap C_{\rm loc}^0(\mathbb{R}^2)$, where
 $$\varphi_0(x)=-\f{1}{4\pi(1-\beta)}\log
 \le(1+\f{\pi}{1-\beta}|x|^{2(1-\beta)}\ri).$$
 Moreover, $c_\epsilon u_\epsilon\rightharpoonup G$ weakly in $W_0^{1,p}(\Omega)$ for any $p$, $1<p<2$,
  strongly in $L^q(\Omega)$ for any $q>1$, and $c_\epsilon u_\epsilon\ra G$ in
 $C^1_{\rm loc}(\overline{\Omega}\setminus\{0\})$, where $G$ is a distributional solution to
 \be\label{g-r}-\Delta G-\alpha G=\delta_0-\sum_{i=1}^m\psi_i(0)\psi_i.\ee Clearly $G$ takes the form
 \be\label{grfunct}G(x)=-\f{1}{2\pi}\log|x|+A_0+{\psi}(x),\ee
where $A_0$ is a constant and ${\psi}\in C^1(\overline{\Omega})$. Since $u_\epsilon\in E_\ell^\perp$, we have
\be\label{G-perp}\int_\Omega G\psi_idx=\lim_{\epsilon\ra 0}\int_\Omega c_\epsilon u_\epsilon \psi_idx=0,\quad
\forall \psi_i\in E_\ell.\ee
Finally we have
\be\label{upper-bd-2}\sup_{u\in W_0^{1,2}(\Omega),\,\|u\|_{1,\alpha}\leq 1}\int_\Omega \f{e^{4\pi(1-\beta)u^2}}{|x|^{2\beta}}dx
=\lim_{\epsilon\ra 0}\int_\Omega\f{e^{4\pi(1-\beta-\epsilon)u_\epsilon^2}}{|x|^{2\beta}}dx\leq
\int_\Omega|x|^{-2\beta}dx+\f{\pi}{1-\beta}e^{1+4\pi(1-\beta)A_0}.\ee
In conclusion, we have proved Theorem \ref{Thm3}. $\hfill\Box$

\subsection{Proof of Theorem \ref{Thm4}}
    In view of (\ref{upper-bd-2}), it suffices to construct a sequence of functions ${\phi}_\epsilon^\ast\in E_\ell^\perp$ with $\|{\phi}_\epsilon^\ast\|_{1,\alpha}=1$
     such that for sufficiently
     small $\epsilon>0$,
     \be\label{gg}
     \int_\Omega \f{e^{4\pi(1-\beta){\phi_\epsilon^\ast}^2}}{|x|^{2\beta}}dx>
     \int_\Omega|x|^{-2\beta}dx+\f{\pi}{1-\beta}e^{1+4\pi(1-\beta)A_0}.
     \ee
     We shall adapt the test functions constructed in the proof of Theorem \ref{Thm2}.
      Let $\phi_\epsilon$ be defined by (\ref{ppp}), where $G$ be defined as in (\ref{g-r}) and (\ref{grfunct}), $R=(-\log\epsilon)
      ^{1/(1-\beta)}$.
      The constants $c$ and $b$ are determined by (\ref{c2-2}) and (\ref{B}) such that
      $\phi_\epsilon$ satisfies
     the following three properties: $(i)$ $\phi_\epsilon\in W_0^{1,2}(\Omega)$; $(ii)$ $\|\phi_\epsilon\|_{1,\alpha}=1$;
     $(iii)$
      there holds
      $$\int_{\Omega}\f{e^{4\pi(1-\beta)\phi_\epsilon^2}}{|x|^{2\beta}}dx\geq\int_\Omega|x|^{-2\beta}dx+\f{\pi}{1-\beta}e^{1+4\pi(1-\beta)A_0}
     +\f{4\pi(1-\beta)}{c^2}
     \int_{\Omega}|x|^{-2\beta}G^2dx+o(\f{1}{c^2}).$$
     Recalling that
     $(\psi_i)_{i=1}^m$ is an orthogonal basis of $E_\ell$ verifying (\ref{orth}), we set
     $$\widetilde{\phi}_\epsilon=\phi_\epsilon-\sum_{i=1}^m(\phi_\epsilon,\psi_i)\psi_i,$$
     where
     $$(\phi_\epsilon,\psi_i)=\int_\Omega \phi_\epsilon \psi_idx.$$
     Obviously $\widetilde{\phi}_\epsilon\in E_\ell^\perp$. A straightforward calculation shows
     \bea
     (\phi_\epsilon,\psi_i)&=&\int_{\mathbb{B}_{R\epsilon}}\le(c+\f{-\f{1}{4\pi(1-\beta)}\log(1+\f{\pi}{1-\beta}\f{|x|^{2(1-\beta)}}
     {\epsilon^{2(1-\beta)}})+b}{c}\ri)
      \psi_idx\nonumber\\[1.2ex]&&+\int_{\mathbb{B}_{2R\epsilon}
     \setminus \mathbb{B}_{R\epsilon}}\f{G-\eta \psi}{c} \psi_idx+
     \int_{\Omega\setminus \mathbb{B}_{2R\epsilon}}\f{G}{c} \psi_idx\nonumber\\
     [1.2ex]&=& o(\f{1}{\log^2\epsilon}).\label{small}
     \eea
     Here we have used (\ref{G-perp}) to derive
     $$\int_{\Omega\setminus \mathbb{B}_{R\epsilon}}\f{G}{c} \psi_idx=-\int_{\mathbb{B}_{R\epsilon}}
     \f{G}{c} \psi_idx=O(\epsilon^2(-\log\epsilon)^{\f{5-\beta}{2(1-\beta)}})
     =o(\f{1}{\log^2\epsilon}).$$
     By (\ref{small}) and property $(ii)$ of ${\phi}_\epsilon$, we have
     \bea
     &&\widetilde{\phi}_\epsilon={\phi}_\epsilon+o(\f{1}{\log^2\epsilon}),\label{333}\\[1.2ex]
     &&\|\widetilde{\phi}_\epsilon\|_{1,\alpha}^2=1+o(\f{1}{\log^2\epsilon}).\label{444}
     \eea
     Combining (\ref{333}), (\ref{444}) and property $(iii)$ of $\phi_\epsilon$, we obtain
     \bna
     \int_{\Omega}\f{e^{4\pi(1-\beta)\f{\widetilde{\phi}_\epsilon^2}{\|\widetilde{\phi}_\epsilon\|_{1,\alpha}^2}}}
     {|x|^{2\beta}}dx&=&\int_{\Omega}\f{e^{4\pi(1-\beta)\phi_\epsilon^2+o(\f{1}{\log\epsilon})}}{|x|^{2\beta}}dx\\
     &\geq&(1+o(\f{1}{\log\epsilon}))\le(\int_\Omega|x|^{-2\beta}dx+\f{\pi}{1-\beta}e^{1+4\pi(1-\beta)A_0}\ri.\\
     &&\quad\quad\quad\le.
     +\f{4\pi(1-\beta)}{c^2}
     \int_{\Omega}|x|^{-2\beta}G^2dx+o(\f{1}{c^2})\ri)\\
     &\geq&\int_\Omega|x|^{-2\beta}dx+\f{\pi}{1-\beta}e^{1+4\pi(1-\beta)A_0}+\f{4\pi(1-\beta)}{c^2}
     \int_{\Omega}\f{G^2}{|x|^{2\beta}}dx+o(\f{1}{c^2}).
     \ena
     Set $\phi_\epsilon^\ast=\widetilde{\phi}_\epsilon/\|\widetilde{\phi}_\epsilon\|_{1,\alpha}$. Since
     $\widetilde{\phi}_\epsilon\in E_\ell^\perp$, we have $\phi_\epsilon^\ast\in E_\ell^\perp$.
     Moreover $\|\phi_\epsilon^\ast\|_{1,\alpha}=1$ and (\ref{gg}) holds. $\hfill\Box$

\bigskip

{\bf Acknowledgements}. Y. Yang is supported by the National Science Foundation of China (Grant No.11171347 and Grant
 No. 11471014). X. Zhu is supported by the National Science Foundation of China (Grant No. 41275063 and Grant
 No. 11401575).

\end{document}